\documentclass[11pt,reqno]{amsart}
\usepackage[T1]{fontenc}
\usepackage[latin9]{inputenc}
\usepackage{geometry}
\usepackage{verbatim}
\usepackage{mathrsfs}
\usepackage{amsthm}
\usepackage{amstext}
\usepackage{amssymb}
\usepackage{graphicx}
\usepackage{esint}
\usepackage{color}
\usepackage{amsfonts}
\usepackage{amsmath} 
\usepackage{bm}           
\usepackage{mathabx}
\usepackage{amsxtra}
\usepackage{stmaryrd}
\usepackage{old-arrows}
\usepackage{float,epsf}
\usepackage[english]{babel}
\usepackage{enumerate}
\usepackage{tikz}
\usepackage{srcltx}
\usepackage{verbatim}
\usepackage[multiple]{footmisc}

\makeatletter

\theoremstyle{plain}
\newtheorem{theorem}{Theorem}[section]

\newcommand{\bu}{{\bf u}}

\def\intave#1{\int_{#1}\hbox{\llap{$\raise2.3pt\hbox{\vrule
height.9pt width7pt}\phantom{\scriptstyle{#1}}\mkern-2mu$}}}

\def\intav#1{\mathchoice
          {\mathop{\vrule width 9pt height 3 pt depth -2.6pt
                  \kern -9pt \intop}\nolimits_{\kern -6pt#1}}%
          {\mathop{\vrule width 5pt height 3 pt depth -2.6pt
                  \kern -6pt \intop}\nolimits_{#1}}%
          {\mathop{\vrule width 5pt height 3 pt depth -2.6pt
                  \kern -6pt \intop}\nolimits_{#1}}%
          {\mathop{\vrule width 5pt height 3 pt depth -2.6pt
                  \kern -6pt \intop}\nolimits_{#1}}}
\def\intav#1{\vint_{#1}}

\allowdisplaybreaks[1]
\allowdisplaybreaks[1]

\newcommand{\divg}{ \mbox{div\,}}

\newcommand{\x}{{\bf x}}

\newcommand{\e}{{\varepsilon}}
\def \epsi {\varepsilon}

\newcommand{\cav}{\mathrm{cav}}

\usepackage{babel}

\usepackage{babel}

\makeatother

\usepackage{babel}

\numberwithin{equation}{section}

\begin{document}
\title[Morawetz's Contributions]{Morawetz's Contributions \\
to the Mathematical Theory of \\Transonic Flows, Shock Waves, and \\ Partial Differential Equations of Mixed Type}
\author{Gui-Qiang G. Chen}
\address{Gui-Qiang G. Chen: Oxford Centre for Nonlinear Partial Differential Equations,
Mathematical Institute, University of Oxford, Oxford, OX2 6GG, UK.
 $\quad$ Email: chengq@maths.ox.ac.uk}

 \begin{abstract}
This article is a survey of Cathleen Morawetz's contributions to the mathematical theory of transonic flows, shock waves,
and partial differential equations of mixed elliptic-hyperbolic type.
The main focus is on Morawetz's fundamental work on the non-existence of continuous transonic flows past profiles,
Morawetz's program regarding the construction of global steady weak transonic flow solutions past profiles via compensated compactness,
and a potential theory for regular and Mach reflection of a shock at a wedge.
The profound impact
of Morawetz's work on recent developments and breakthroughs
in these research directions and related areas in pure and applied mathematics are also discussed.
\end{abstract}
\maketitle

\section{Introduction}
It is impossible
to review all of Cathleen Morawetz's
paramount contributions to pure and applied mathematics and to fully assess their impact on twentieth
century mathematics and the mathematical community in general.
In this article, we focus on Morawetz's deep and influential work on the analysis of partial differential equations (PDEs)
of mixed elliptic-hyperbolic type,
most notably in the mathematical theory of transonic flows and shock waves.
We also discuss the profound impact
of Morawetz's work on some
recent developments and breakthroughs
in these research directions and related areas in pure and applied mathematics.

\vspace{2pt}
Morawetz's early work on transonic flows
has not only provided a new understanding of mixed-type PDEs,
but has also led to new methods of efficient
aircraft design.
Morawetz's program for constructing global steady weak transonic flow solutions
past profiles has been a source of motivation for numerous recent developments
in the analysis of nonlinear PDEs
of mixed type and related mixed-type problems through weak convergence methods.
Furthermore, her work on the potential theory for regular and Mach reflection of a shock at a wedge
(now known as the von Neumann problem)
has been an inspiration
for the recent complete solution of the von Neumann conjectures regarding global
shock regular reflection-diffraction configurations, all the way up to the detachment angle of the wedge.

\vspace{2pt}
As a graduate student, I learned a great deal from Cathleen's papers [20--25],
which were a true inspiration to me.
My academic journey took a significant turn when I joined the Courant Institute of Mathematical Sciences (New York University)
as a postdoctoral fellow under the direction of Peter Lax. During this time, I had the
extraordinary opportunity to learn directly from Cathleen about
the challenging and fundamental research field that had, until that time, remained largely unexplored.
I was immensely grateful to Cathleen for dedicating countless hours to discuss and analyze with me
a long list of open problems in this field. Her insights were both illuminating and prolific,
and I learned immensely from her during my these years at Courant.
Making substantial progress on some of these longstanding open problems, however,
was a journey that spanned over 10 years, on and off.
Indeed, this field has proven to be truly
challenging.
As a result, I experienced great joy when I had the honor of presenting our first
solution of the von
Neumann problem in \cite{CF2010} to Cathleen
during my lecture at the Conference on Nonlinear Phenomena in Mathematical Physics,
dedicated to her
on the occasion of her 85th birthday, held at the Fields Institute in Toronto, Canada,
from the 18th to 20th of September 2008.

\section{Background}

The work of Cathleen Morawetz on {\it transonic flows}
and the underlying PDEs of mixed elliptic-hyperbolic type
spanned her career.
Let us first describe some of the background regarding {\it transonic flows}
in aerodynamics.

\vspace{2pt}
A fundamental question in aerodynamics is whether it is realistic
for an aircraft to fly at a relatively high speed, with respect
to the speed of sound in the surrounding air,
with both relatively low economical and
environmental costs.
It is known that, at relatively low speeds -- the subsonic
range, the wing can be {\it sailed} by designing it
to {\it obtain as much of a free ride} as possible from the wind.
At very high speeds -- the supersonic range, {\it rocket propulsion} is needed to
overcome the drag produced by shocks that invariably
form (the sonic boom).
The purpose of studying transonic
flows is to find a compromise that allows for {\it sailing}
efficiently {\it near the speed of sound}; this is a critical speed at which
aerodynamic challenges emerge due to shock formation and increased drag ({\it cf}. \cite{Mo82}).

\vspace{2pt}
As shown in Fig. 1, shocks (depicted in bold black) begin to
appear on a wing when the wing speed is below, but near, the speed of sound.
As the wing speed increases from
subsonic (Mach number $M < 1$) to supersonic
(blue region,  $M > 1$), some supersonic shock
(depicted in bold black) appears over the wing already at the Mach number $M=0.77$.
As the wing speed increases from subsonic ($M<1$) to supersonic ($M>1$),
additional shocks and transonic regions are formed.

\begin{figure}
\vspace{-8pt}
\begin{minipage}{0.45\textwidth}
\centering
\includegraphics[height=1.50in,width=2.40in]{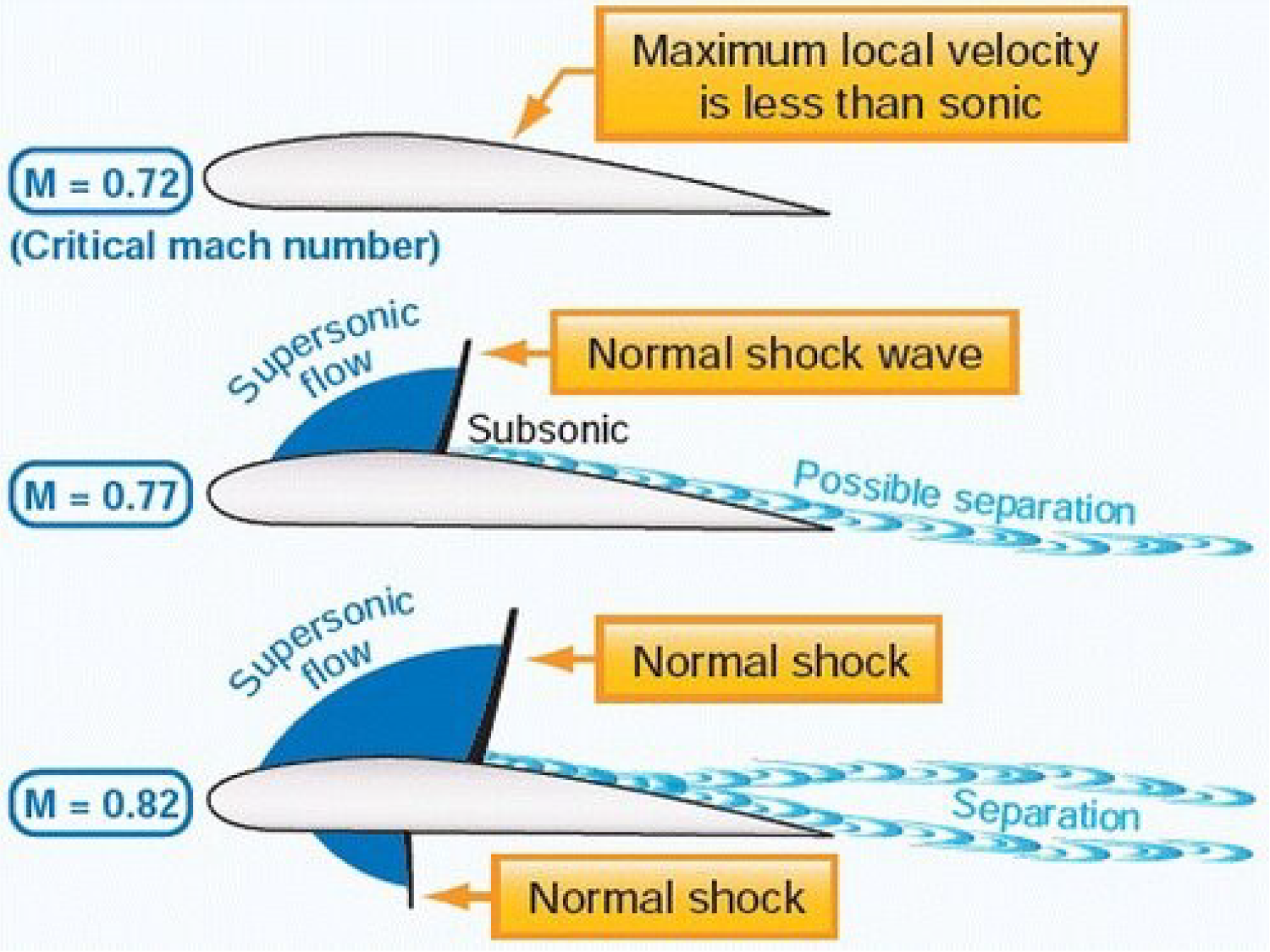}
\end{minipage}
\begin{minipage}{0.51\textwidth}
\centering
\includegraphics[height=1.50in,width=2.65in]{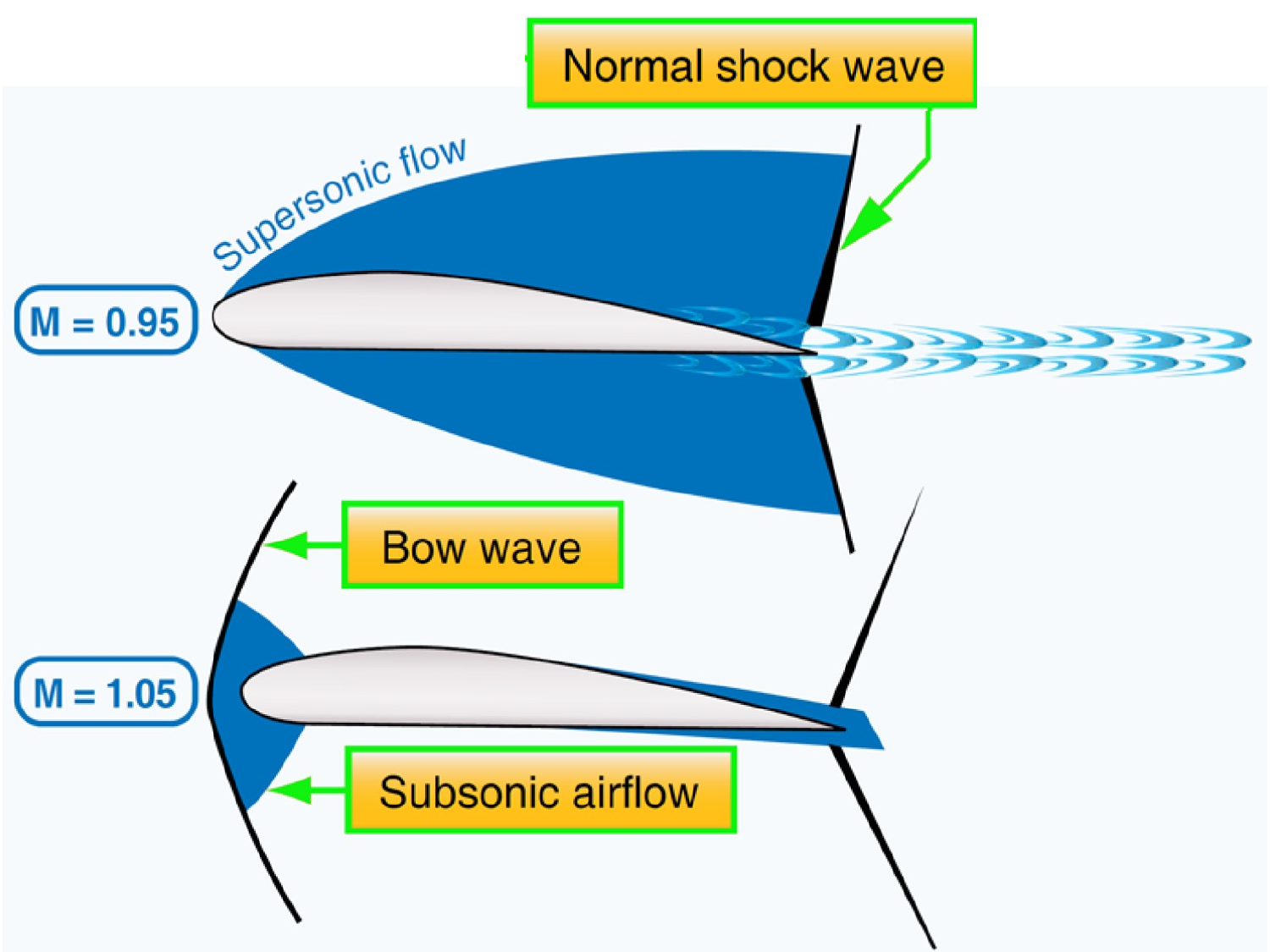}
\end{minipage}
\caption[]{\small
Transonic flow patterns around a wing at and above the critical Mach numbers.
As the wing speed increases from subsonic (Mach number $M<1$) to supersonic (Mach number $M>1$),
some supersonic shock (depicted in bold black) appears over the wing already at the Mach number $M=0.77$.}
\end{figure}

\vspace{2pt}
The two-dimensional irrotational, stationary, compressible and isentropic flow of air around
a profile $\mathcal{P}$ is governed by the Euler equations of the conservation law of
mass and Bernoulli's law for the velocity potential $\varphi(x,y)$
and the density $\rho$ of the fluid ({\it cf.} \cite{Bers,CF18,CFr,Da,GM}):
\begin{equation}\label{potential-eq}
\divg\big(\rho \nabla\varphi\big)=0, \qquad
{1\over 2}|\nabla \varphi|^2+\frac{\rho^{\gamma-1}}{\gamma-1}=B.
\end{equation}
Here $\rho$ is the density, $\varphi$ is the velocity potential ({\it i.e.}, $\nabla\varphi=(\varphi_x, \varphi_y)$ is the velocity),
$\gamma>1$ is the adiabatic exponent for the ideal gas ($\gamma\approx 1.4$ for
the usual gas), and $B>0$ is the Bernoulli constant.
Eq. \eqref{potential-eq} can be formulated as the following system for the velocity $\bu=(u,v)=\nabla \varphi$:
\begin{equation}\label{potential-system-1}
\begin{cases}
v_x-u_y =0,\\[1mm]
(\rho_B(|\bu|) u)_x +(\rho_B(|\bu|) v)_y=0,
\end{cases}
\end{equation}
with
\begin{equation}\label{potential-system-2}
\rho_B(q)=\big(B-\frac{1}{2}q^2\big)^{\frac{1}{\gamma-1}}.
\end{equation}

\vspace{2pt}
System \eqref{potential-eq} for the steady
velocity potential $\varphi$
can be rewritten as
\begin{equation}\label{steady-potential-eq}
\divg\big(\rho_B(|\nabla\varphi|)\nabla\varphi\big)=0,
\end{equation}
or, equivalently, as
\begin{equation}\label{18-b}
a_{11}\varphi_{xx} +2a_{12}\varphi_{xy}
+ a_{22}\varphi_{yy}=0,
\end{equation}
with $a_{11}=c^2-\varphi_x^2$, $a_{12}=-\varphi_x\varphi_y$, and $a_{22}=c^2-\varphi_y^2$,
where $c>0$ is the local speed of sound defined by
$$
c^2=\frac{{\rm d} p}{{\rm d}\rho}
=\rho^{\gamma-1}
=(\gamma-1)\rho_B^{\gamma-1}=(\gamma-1)\big(B-\frac{1}{2}|\nabla \varphi|^2\big).
$$
The adiabatic pressure-density relation for the air is
$p=p(\rho)=\rho^\gamma/\gamma$ (after scaling) with $\gamma\approx 1.4$.

\vspace{2pt}
Eq.\,\eqref{steady-potential-eq} is a {\it nonlinear conservation law of
mixed elliptic-hyperbolic type} for the velocity potential $\varphi$ ({\it cf}. \cite{Bet,Chen2023,CF22}); that is, it is
\begin{itemize}
\item  strictly {\it elliptic} (subsonic), {\it i.e.}, $a_{12}^2-a_{11}a_{22}<0$,
if
$|\nabla\varphi|<c_*:=\sqrt{\frac{2(\gamma-1)B}{\gamma+1}}$;

\smallskip
\item strictly {\it hyperbolic} (supersonic) if
$|\nabla\varphi|>c_*$.
\end{itemize}
The transition boundary
here is $|\nabla\varphi|=c_*$ (sonic),
a degenerate set of \eqref{steady-potential-eq},
which is {\it a priori} unknown,
as it is determined by the gradient of the solution itself.

\vspace{2pt}
The natural boundary condition for the airfoil problem is
\begin{equation}\label{18}
\nabla\varphi \cdot \mathbf{n} = 0 \qquad \mbox{on $\partial\mathcal{P}$},
\end{equation}
where $\mathbf{n}$ is the unit inner normal on the airfoil boundary  $\partial\mathcal{P}$.
The nature of the flow is determined by the local Mach
number $M = \frac{q}{c}$, where $q=|\nabla\varphi|$ is the flow speed.

\vspace{2pt}
A transonic flow occurs when the flow
involves both subsonic and supersonic regions, and Eq. \eqref{steady-potential-eq}
is of mixed elliptic-hyperbolic type.
In this context,
supersonic regions are characterized by the presence of shocks that result from drastic changes in the air density and pressure, due to compressibility.
These pressure changes propagate at supersonic speeds,
giving rise to a shock wave, typically with a small but finite thickness.
In Fig. 1, the shock wave is depicted in {bold black}.
The velocity field $\nabla\varphi$, governed by \eqref{steady-potential-eq}, experiences
jump discontinuities as the flow crosses the shock wave.
These discontinuities serve as an indicator of the presence of shocks.
The mathematical description of
shocks involves an analysis of entropy effects ({\it cf}. \cite{CF18,CFr,Da,Lax}).
The corresponding entropy condition states, mathematically, that {\it  the density function
increases across the shock in the flow direction} ({\it cf}. \cite{Chen2023,CF18,CF22,CFr}).

\section{The Non-Existence of Continuous Transonic Flows past Profiles}

Between the 1930s and the 1950s, there was a long debate among leading scientists, including G. I. Taylor and A. Busemann,
regarding transonic flows around a given airfoil.
The central questions were:
\begin{itemize}
\item Do transonic flows about a given airfoil always, never, or sometimes produce shocks?
\item Is it possible to design a viable airfoil capable of shock-free
flight at a range of transonic speeds?
\end{itemize}
There was no
definitive, satisfactory answer to these questions until Morawetz's work in the 1950s.
In a series of papers \cite{Mo56,Mo57,Mo58},
Morawetz provided a mathematically definite answer to the questions by
proving that shock-free transonic flows are unstable with
respect to arbitrarily small perturbations in the shape
of the profile. More precisely, this can be stated as follows:

\smallskip
\begin{theorem}[Morawetz \cite{Mo56,Mo57,Mo58}]
{Let $\varphi$ be a transonic solution of \eqref{steady-potential-eq}--\eqref{18} with
a continuous velocity field $\nabla\varphi$ and a fixed speed $q_\infty$ at infinity
about a symmetric profile $\mathcal{P}$ $(${\it cf.} {\rm  Fig. 1}$)$. For an arbitrary
perturbation $\tilde{\mathcal{P}}$ of $\mathcal{P}$ along an arc inside the supersonic
region attached to the profile that contains the point of
maximum speed in the flow, there is no continuous velocity field $\nabla\tilde{\varphi}$
solving the corresponding problem \eqref{steady-potential-eq}--\eqref{18} with $\tilde{P}$.}
\end{theorem}

\smallskip
This indicates that any arbitrary perturbation of the airfoil inside the blue supersonic
region creates shocks.
In simpler terms,
even if a viable profile capable of a shock-free transonic
flow can be designed, any imperfection in its construction leads to the
formation of shocks at the intended speed.

\vspace{2pt}
Morawetz's proof involved ingenious new estimates for the solutions of
nonlinear PDEs of mixed elliptic-hyperbolic type.
The proof made two major advances:
First,
the correct boundary value problem was formulated,
satisfying the perturbation of the velocity potential
in the hodograph
plane. In this plane, a hodograph transformation linearized
Eq. \eqref{steady-potential-eq}
and mapped the known profile exterior into an
unknown domain.
Second, by developing carefully tailored integral
identities, a uniqueness theorem was proven for regular
solutions of the transformed PDE, with data prescribed
on only a proper subset of the transformed boundary
profile. This theorem states that the transformed problem is
predetermined, and that  no regular solutions exist.
This result was further extended to include fixed profiles but finite
perturbations in $q_\infty$, as well as to non-symmetric profiles (see \cite{Cook}).

\vspace{2pt}
Morawetz's results catalyzed a significant change in the views
of engineers,
compelling them
to re-calibrate wing design to minimize
the shock strength over a useful range of transonic speeds.
The work of H.~H.
Pearly and later R.~R. Whitecomb in the 1960s on supercritical airfoils
further underscored
the impact of
Morawetz's findings on
transonic airfoil design.
In the 1970s, this research direction experienced a surge in growth
as part of the field of computational fluid dynamics. Key milestones
included the type-dependent difference scheme proposed by E.~M. Merman and
J.~D. Cole (1971), the complex characteristic method proposed by P.
Paramedian and D. Koran (1971), and the rotated difference
scheme introduced by A. Jameson (1974).
These advances significantly contributed to the accurate calculation
of steady transonic flows and the development of codes for transonic airfoil design ({\it cf}. \cite{Jameson}).

\vspace{2pt}
Cathleen Morawetz's work on transonic flows not only
transformed the field of PDEs of mixed elliptic-hyperbolic type, but also
served as a compelling example of mathematics coming to the rescue with regard to
real-world problems.
At the time, when many engineers and applied scientists were deeply sceptical
about the role of mathematics
in terms of real-world applications,
Morawetz's work demonstrated the true usefulness of the discipline.

\vspace{2pt}
In the most recent two decades, Morawetz's work has served as a source of inspiration for mathematicians and has led
to numerous
significant developments in the
analysis of steady transonic flows and free boundary problems for the steady compressible
Euler equations and other nonlinear PDEs of mixed elliptic-hyperbolic type.
These developments encompass a wide range of transonic flow scenarios,
including those around wedges and conical bodies,
transonic nozzle flows (including the de Laval nozzle flow), and other related steady transonic flow problems.
For more details on these developments,
we refer the reader to \cite{Chen2023,CF18,CF22} and the references provided therein.

\smallskip
\section{Morawetz's Program for the Construction of Global Weak Transonic Flows past Profiles via Compensated Compactness}

With the complete solution regarding the non-existence of continuous transonic flows
past profiles, {\it i.e.},
the {\it exceptional nature} of shock-free transonic flows, Morawetz
turned to the next fundamental problems:
\begin{itemize}
\item {Can robust existence theorems for weak shock solutions be established?}
\item {Can a weak shock be contracted to a sonic point on the profile?}
\end{itemize}
The first problem, now known as the Morawetz problem,
received support from the work of Garabedian-Korn
in 1971, which demonstrated that small perturbations of
continuous flows can result in only weak shocks in the case of potential flow.
The second
question was inspired by the thinking of K.~G. Guderley in
the 1950s.

\vspace{2pt}
Inspired by the difference method of A. Jameson (1974),
Morawetz introduced an artificial viscosity parameter
into the nonlinear potential equation in \cite{Mo85,Mo94b}.
This viscosity method involved replacing
the Bernoulli law, which previously related density $\rho$ to the gradient of the velocity potential $\nabla\varphi$
as
$
\rho=\rho_B(|\nabla\varphi|),
$
with a first-order PDE that retards the density $\rho$.
Morawetz presented an ambitious program aimed at proving the existence of global weak solutions of the problem.

This program involved embedding the problem
within an assumed viscous framework, wherein
the following compensated compactness
framework would be satisfied:

\smallskip
\begin{theorem}[Morawetz \cite{Mo85,Mo94b}]\label{thm4.1}
Let
$\{\mathbf{u}^\e:=\nabla\varphi^\e\}_{\e>0} \subset L^\infty(\Omega)$
be a sequence of approximate solutions
of the Morawetz problem for \eqref{potential-system-1}--\eqref{potential-system-2} in the domain $\Omega$
with the following uniform bounds{\rm :}
\begin{itemize}
\item[(i)] There exist $q_*, q^*\in (0, q_{\cav})$
independent of $\epsi>0$
such that
\begin{equation}\label{5.2a}
    0 < q_* \leq |\bu^\varepsilon(\x)| \leq q^*
    \qquad \mbox{for all $\x\in \Omega$},
\end{equation}
where $q_{\cav}$ is the maximum speed so that $\rho_B(q_{\cav})=0$.

\item[(ii)] There are $\theta_*,\theta^*\in (-\infty, \infty)$ such that the corresponding velocity angle sequence
$\theta^\e(\x)$ satisfies
\begin{equation}\label{5.3a}
    \theta_*\le \theta^\varepsilon(\x)\le \theta^*
    \qquad \mbox{for all $\x\in \Omega$}.
\end{equation}

\item[(iii)]
The corresponding entropy dissipation measure sequence
\begin{equation}\label{5.3b}
{\rm div}_\mathbf{x}(Q_1, Q_2)(\bu^\epsi)
\quad\mbox{is compact in $H^{-1}_{\rm loc}(\Omega)$},
\end{equation}
where $(Q_1, Q_2)$ is any $C^2$ entropy pair of system  \eqref{potential-system-1}--\eqref{potential-system-2}, and
$H^{-1}(\Omega)$ is the dual space of the Sobolev space $W^{1,2}(\Omega):=H^1(\Omega)$.
\end{itemize}
Then there exist a subsequence $($still denoted$)$ $\bu^\epsi$ and a function $\bu\in L^\infty(\Omega)$ such that
\begin{equation}\label{eqn-convergence-approx}
\bu^\epsi\,\to \,\bu
\qquad\quad \mbox{pointwise {\it a.e.}, and in $L^p$ for any $p\in [1,\infty)$}.
\end{equation}
\end{theorem}

\smallskip
With this framework, Morawetz demonstrated that the sequence of solutions of the viscous problem,
which remains uniformly away from both stagnation ($|\bu|=0$) and cavitation ($\rho=0$)
with uniformly bounded velocity angles,
converges subsequentially to an entropy solution
of the transonic flow problem; see also Gamba-Morawetz \cite{GaM}.
Some possible extensions of this framework were also discussed in \cite{Mo85,Mo94b}.

\vspace{2pt}
An alternative vanishing viscosity method,
building upon
Morawetz's pioneering work, was developed in \cite{CSW}.
This method was designed for adiabatic constant $\gamma\in (1,3)$
and ensures a family of invariant
regions for the corresponding viscous problem. This implies an upper bound,
uniformly away from cavitation, for the viscous approximate velocity fields.
In other words,
the condition in \eqref{5.2a},
\begin{equation}\label{5.2b}
    |\bu^\epsi (\x)| \leq q^*
    \qquad \mbox{for all $\x\in \Omega$},
\end{equation}
can be verified rigorously  for the viscous approximate solutions for $\gamma\in (1,3)$.
This method involves the construction of mathematical entropy pairs through the Loewner-Morawetz relation via
entropy generators governed by a generalized Tricomi equation of mixed elliptic-hyperbolic type.
The corresponding entropy dissipation measures are analyzed
to ensure that the viscous solutions satisfy the compactness
framework (Theorem \ref{thm4.1}).
Consequently, the compensated compactness framework implies that
a sequence of solutions to the viscous problem,
staying uniformly away from stagnation with uniformly bounded velocity angles,
converges to an entropy solution
of the inviscid transonic flow problem.

\vspace{2pt}
On the other hand, for the case $\gamma\ge 3$, cavitation does occur.
In particular, for $\gamma=3$ (similar to the case $\gamma=\frac{5}{3}$ for the time-dependent isentropic case),
a family of invariant regions for the corresponding viscous problem has been identified in the most recent preprint \cite{CST}.
This identification implies a lower bound
uniformly away from stagnation, but it allows for the occurrence of cavitation, in
the viscous approximate velocity fields.
As a result, the first complete existence theorem for system \eqref{potential-system-1}--\eqref{potential-system-2}
has been established, without requiring {\it a priori} estimate assumptions.

\section{Potential Theory for Regular and Mach Reflection of a Shock at a Wedge}

When a planar shock
separating two constant states, (0) and (1),
with constant velocities
and densities $\rho_0<\rho_1$
(state (0) is ahead or to the right of the shock, and state
(1) is behind the shock), travels in the flow direction
and impinges upon a symmetric wedge
with a half-wedge angle $\theta_{\rm w}$
head-on at time $t = 0$,
a reflection-diffraction process occurs as time progresses ($t > 0$).
A fundamental question arises regarding the nature of the wave patterns
formed in
shock reflection-diffraction configurations around the wedge.
The complexity of these
configurations was initially reported
by Ernst Mach (1878),
who observed two patterns of
shock reflection-diffraction configurations:
Regular reflection (characterized by a two-shock configuration)
and Mach reflection (characterized by a three-shock/one-vortex-sheet configuration).
These configurations are illustrated in Fig. 2, below\footnote{ M. Van Dyke:
{\it An Album of Fluid Motion}, The Parabolic Press: Stanford, 1982.}.
The issue remained largely unexplored until the 1940s when John von Neumann \cite{vN1,vonN},
along with other mathematical and experimental
scientists (see, for example,  \cite{BD,CF18,CFr,FS,GM}
and the references cited therein), initiated
extensive research into all aspects of shock reflection-diffraction phenomena.

\begin{figure}[h]\label{cog}
  \centering
  \includegraphics[width=0.18\textwidth]{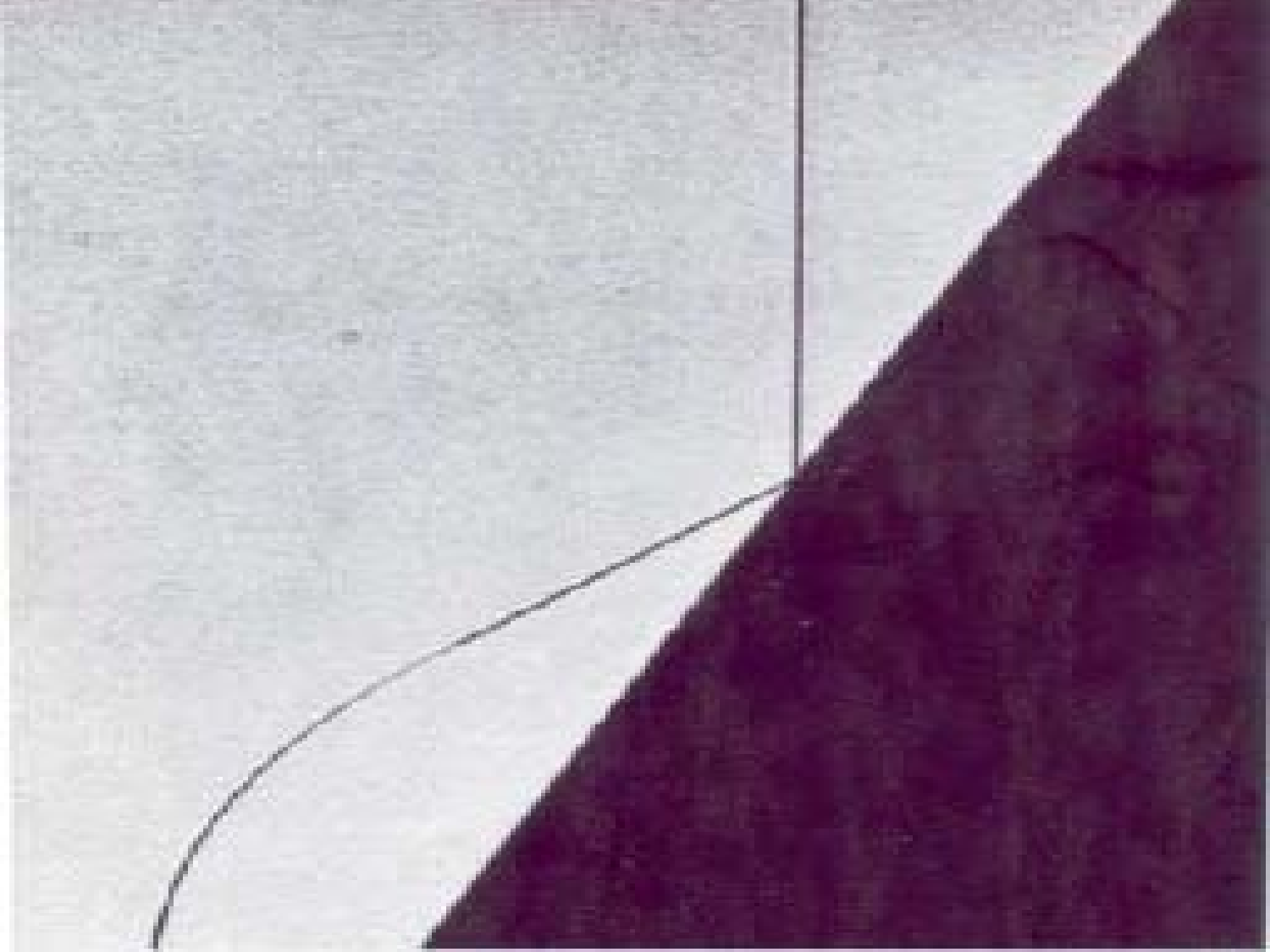} \qquad
  \includegraphics[width=0.18\textwidth]{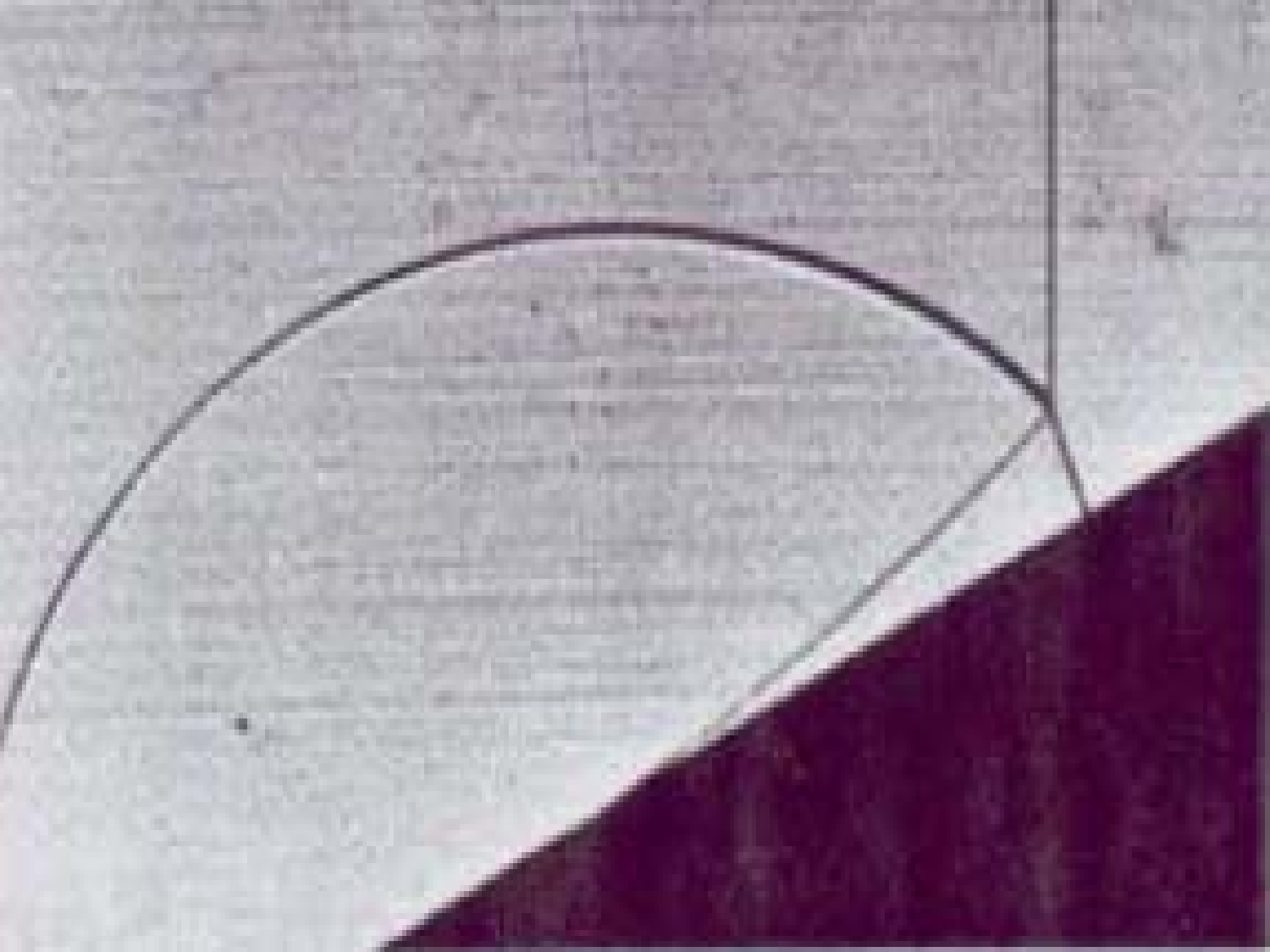}\qquad
  \includegraphics[width=0.18\textwidth]{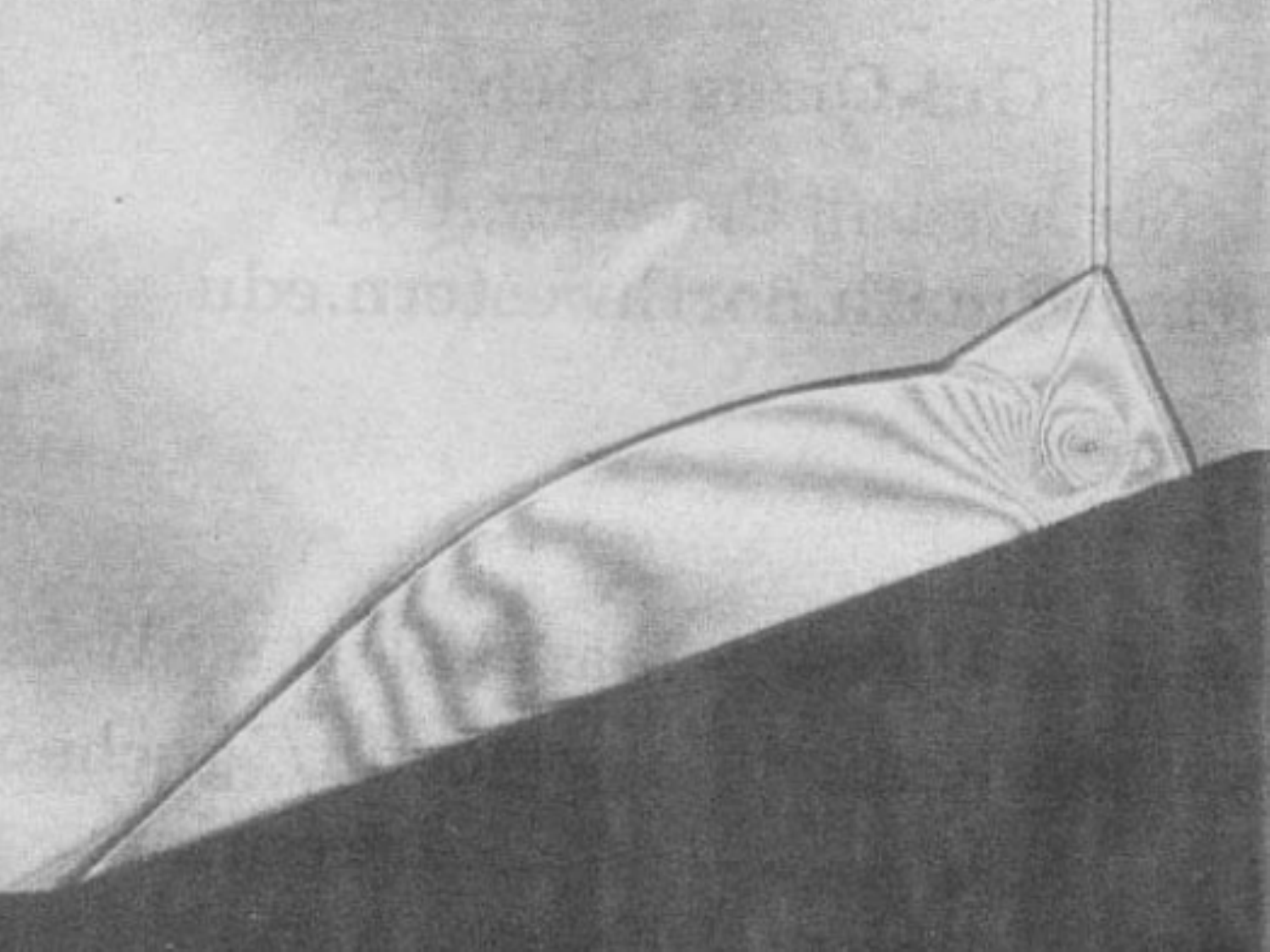}
  \caption[]{\small Three patterns of shock reflection-diffraction configurations.}
\end{figure}

\vspace{2pt}
In fact, the situation is far more intricate than
what Mach initially observed.
Shock reflection can be further
divided into more specific sub-patterns, and numerous other patterns of
shock reflection-diffraction configurations exist. These include
supersonic regular reflection, subsonic regular reflection,
attached regular reflection, double Mach reflection,
von Neumann reflection, and Guderley reflection.
For a comprehensive exploration of these patterns, we refer to
\cite{BD,Chen2023,CF18,CFr,GM}
and the references cited therein (see also Figs. 2--4).
The fundamental scientific issues related to shock reflection-diffraction configurations encompass
the following:
\begin{itemize}
\item[(i)] Understanding the structures of these configurations;

\item[(ii)] Determining the transition criteria between the different patterns;

\item[(iii)] Investigating the dependence of these patterns on physical parameters such as
the incident-shock-wave Mach number ({\it i.e.}, the strength of the incident shock),
the wedge angle $\theta_{\rm w}$,
and the adiabatic exponent $\gamma>1$.
\end{itemize}

\begin{figure}
\centering
\begin{minipage}{0.470\textwidth}
\centering
\includegraphics[height=1.5in,width=2.1in]{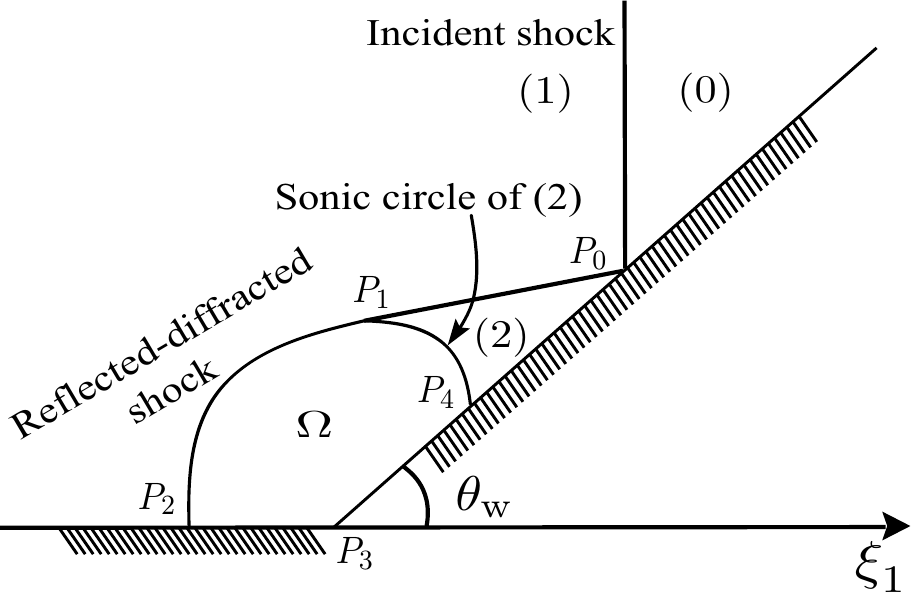}
\caption{\scriptsize \,\, Supersonic regular reflection-diffraction configuration \cite{CF18}.}
\vspace{-2pt}
\label{fig:RegularReflection}
\end{minipage}
\,\,
\begin{minipage}{0.470\textwidth}
\centering
\includegraphics[height=1.5in,width=2.1in]{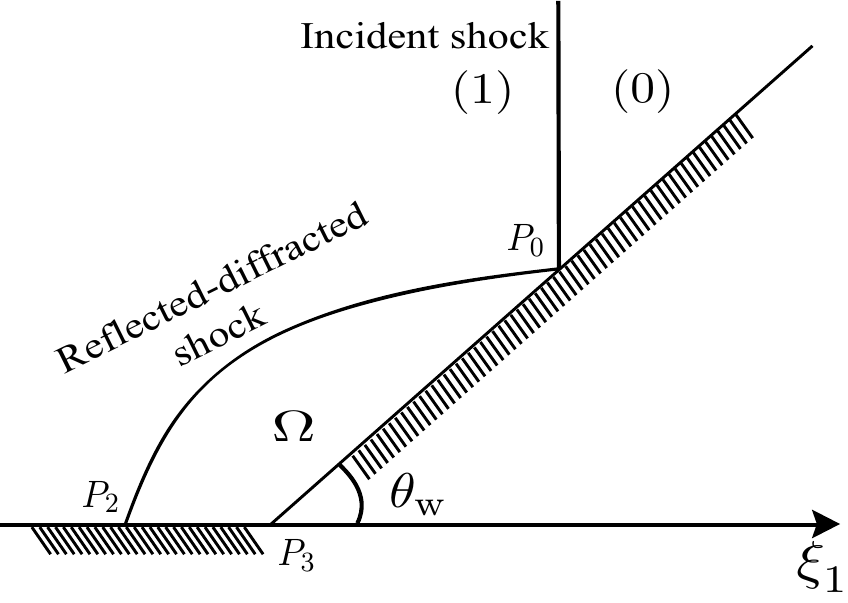}
\caption{\scriptsize \,\, Subsonic regular reflection-diffraction configuration \cite{CF18}.}
\vspace{-2pt}
\label{fig:SubsonicRegularReflection}
\end{minipage}
\vspace{-7pt}
\end{figure}

\vspace{2pt}
\noindent
In particular, several transition criteria between different
patterns of shock reflection-diffraction configurations have been proposed.
Notably, there are two significant conjectures: the {\it sonic conjecture} and the {\it detachment conjecture},
introduced by von Neumann
\cite{vN1} (these are also discussed in \cite{BD,Chen2023,CF18}).

\vspace{2pt}
In her pioneering work \cite{Mo94}, Morawetz investigated
the nature of shock reflection-diffraction patterns for weak incident shocks with strength
$M_1$ that represents the Mach number of state (1) behind the incident shock,
and for small half-wedge angles $\theta_{\rm w}$.
Her exploration involved the use of potential theory, several scaling techniques, the study of
mixed equations, and matching asymptotics for different scalings.
The self-similar equations encountered in the analysis are of
mixed elliptic-hyperbolic type.
The investigation included linearization to obtain a linear mixed flow valid away from a sonic curve.
Near the sonic
curve, a shock solution was constructed by using a different scaling, except near the zone of interaction between
the incident shock and the wall, where a special scaling was employed.
The analysis specifically addressed polytropic gases with an adiabatic exponent $\gamma>1$.
The parameter $\beta=\frac{\theta_{\rm w}^2}{(\gamma + 1)M_1}$
spans a range from $0$ to $\infty$, and its values have distinct implications:

\begin{itemize}
\item When $\beta>2$, regular reflection
(either weak or strong)
is possible, and the entire pattern
is reconstructed to the
lowest order in shock strength;

\item When $\beta<\frac{1}{2}$,  Mach reflection occurs, and the flow behind the reflection
becomes subsonic which, in principle, can be constructed (with an open elliptic problem) and matched;

\item For values of $\beta$ between $\frac{1}{2}$ and $2$, or an even larger value of $\beta$, the flow behind
a Mach reflection may become
transonic. Further investigation is required to determine the nature of this transonic flow.
\end{itemize}

The fundamental pattern of shock reflection includes an almost semi-circular shock.
In the case of regular reflection, this shock originates from the reflection point on the wedge.
In the case of Mach reflection, it is matched with a local interaction flow.
Given their nature, choosing the least entropy generation, the weak regular reflection occurs
for a sufficiently large value of $\beta$ (which settled the von Neumann paradox).
An accumulation point of vorticity is observed on the
wedge, located above the leading point.

\vspace{2pt}
Inspired by Morawetz's work, significant progress has been made on the von Neumann problem for shock-reflection-diffraction.
In particular, several new ideas, approaches, and techniques have been developed to solve fundamental open problems
involving transonic shocks and related free boundary problems for nonlinear PDEs of mixed hyperbolic-elliptic
type; see \cite{Chen2023,CF2010,CF18,CF22} and the references cited therein.
In \cite{CF2010}, a new approach was first introduced and corresponding techniques
were developed to solve the global problem of shock reflection-diffraction by large-angle wedges.
This development eventually led to the first complete and rigorous proof of both von Neumann's sonic and detachment conjectures
in Chen-Feldman \cite{CF18}. Consequently,
this provided a rigorous proof of the existence, uniqueness, stability, and optimal regularity of global solutions
of the von Neumann problem for shock reflection-diffraction, all the way up to the detachment angle,
for the two-dimensional Euler equations for potential flow.
The approaches and techniques developed in these works have also been applied to solve other transonic shock
problems and related problems involving similar difficulties for nonlinear PDEs of mixed type, which arise in
fluid dynamics and other fields (see, for example, the references provided in \cite{Chen2023,CF2010,CF18,CF22}).
Notably, these developments include the first complete
solution to two other longstanding open problems: the Prandtl conjecture for the Prandtl-Meyer configuration
for supersonic flow onto a solid ramp up to the detachment angle (see \cite{BCF-14}) and the global solutions
of the two-dimensional Riemann problem with four-shock interactions for the Euler
equations (see \cite{CCHLW}).

\medskip
\section{Concluding Remarks}
From the brief discussion above,
we have seen that Cathleen Morawetz has made truly seminal contributions to the mathematical theory of transonic flows,
shock waves, and PDEs of mixed elliptic-hyperbolic type.
Furthermore, Morawetz also made fundamental contributions to many other fields, including functional inequalities and scattering theory.
Her work has resulted in numerous significant developments and breakthroughs in various research directions
and related areas in pure and applied mathematics.
These new developments underscore the enduring importance of pioneering work in mathematics and its lasting impact
on contemporary research, further highlighting the significance of Cathleen Morawetz's contributions to these fields.
It is hard to imagine a world without Morawetz's profound contributions to transonic flow theory, shock wave theory,
PDEs of mixed type, functional inequalities, and scattering theory.

\vspace{2pt}
Cathleen Morawetz's role as a leader, mentor, and inspirational figure in the mathematical community has been
equally significant throughout her ninety-four years of life.
Her dedication to excellence and her generosity in sharing her knowledge and insights have left an indelible mark
on generations of mathematicians.
Her legacy will undoubtedly continue to inspire and shape the future of mathematical research for many decades
to come.

\medskip
\bigskip
\section*{About the Author}
Gui-Qiang G. Chen is currently the Statutory Professor in the Analysis of Partial Differential Equations,
the Director of the Oxford Centre for Nonlinear Partial Differential Equations (OxPDE) at the Mathematical Institute,
and a Professorial Fellow of Keble College at the University of Oxford.
His main research areas are in partial differential equations (PDEs), nonlinear analysis, and their applications to other areas of mathematics and science.
His recent research interests include nonlinear hyperbolic systems of conservation laws, nonlinear PDEs of mixed type, nonlinear waves, free boundary problems,
geometric analysis, and stochastic PDEs. His research interests also include measure-theoretical analysis, weak convergence methods,
entropy analysis, statistical physics, and numerical analysis. He is a member of the Academia Europaea, a fellow of the European Academy of Sciences, a fellow of the American Mathematical Society,
and a fellow of the Society of Industrial and Applied Mathematics.

\end{document}